\input amstex
\documentstyle{amsppt}
\magnification=\magstep1

\topmatter
\title{On the space of morphisms between \'etale groupoids } \endtitle  \author  Andre Haefliger \endauthor

\abstract{Given two \'etale groupoids $\Cal G$ and $\Cal G'$, we consider the set of pointed morphisms from $\Cal G$ to $\Cal G'$. Under suitable hypothesis we introduce on this set a structure of Banach manifold  which can be considered as the space of objects of an \'etale groupoid whose space of orbits is the space of morphisms from $\Cal G$ to $\Cal G'$.} \endabstract

\endtopmatter

 This is a drawing up of a talk given at the IHP Paris in february 2007, in the framework of the trimester "Groupoids and Stacks". In this revised version, we correct the paragraph IV.4 and  give more details on the proofs.

 Let $\Cal G$ and $\Cal G'$ be two \'etale groupoids. 
The aim of this talk is to show  that, under suitable conditions, the set of continuous morphisms from $\Cal G$ to $\Cal G'$ can be considered as the space of orbits of  an \'etale groupoid. This was done and used in the paper {\it Closed geodesics on orbifolds} \cite 4 by K. Guruprasad and A. Haefliger  in the particular case where $\Cal G$ is the circle (considered as a trivial groupoid) and $\Cal G'$ an orbifold. The same problem has been also considered independentlly by Weimin Chen in \cite 2. See also the paper \cite 8 by E. Lupercio and B. Uribe.
Our basic references are \cite 1 and \cite 5.
\subheading{Plan of the talk}
\medskip
I.  \'Etale groupoids (homomorphisms, equivalence, localization, developability)

II. Morphisms and cocycles 

III. Pointed morphisms

IV. Topology on the set of pointed morphisms

V. Selfequivalences and extensions.

\heading{I. \'etale groupoids}  \endheading
 
\subheading{1. Definition} An \'etale groupoid $\Cal G$ with space of objects $T$ is a topological groupoid such that the source and target projections $s,t : \Cal G \to T$ are \'etale maps (local homeomorphisms). For $x \in T$, we note $1_x \in \Cal G$ the corresponding unit element.  An arrow $g \in \Cal G$ with source $s(g)=x$ and target $t(g) =y$ is pictured as an arrow $y \leftarrow x$.  Accordingly, the space  of composable arrows, namely the subspace of $\Cal G \times \Cal G$ formed by the pairs $(g,g')$ with $s(g) = t(g')$, is noted $\Cal G \times_T \Cal G$.
We shall  assume throughout that the space of objects $T$ is arcwise locally connected and locally simply connected. The groupoid $\Cal G$ is connected iff its space of orbits $\Cal G \backslash T$ is connected.

\subheading{2. Homomorphisms and equivalences} If $\Cal G'$ is an \'etale groupoid with space of objects $T'$, a homomorphism $\phi:\Cal G \to \Cal G'$ is a continuous functor. It induces a continuous map $\phi_0:T \to T'$, and passing to the quotient,  a continuous map $\Cal G \backslash  T\to \Cal G'\backslash T'$ of spaces of orbits. Two homomorphisms $\phi$ and $\phi'$ from $\Cal G$ to $\Cal G'$ are equivalent if they are related by a natural transformation, i.e a continuous map $h:T \to \Cal G'$ such that the following is defined and satisfied for each $g\in \Cal G$:$$\phi'(g) = h(t(g))\phi(g)h^{-1}(s(g)).$$

We say that $\phi$ is an equivalence if the induced map on the spaces of orbits is bijective and if $\phi$ is locally an isomorphism. This means that each point has an open neighbourhood $U$ such that the restriction of $\phi_0$ to  $U$ is a homeomorphism onto an open set $U'$ of $T'$ and that $\phi$ restricted to $U$ is an isomorphism from the restriction of $\Cal G$ to $U$ to the restriction  of $\Cal G'$ to $U'$. This generates an equivalence relation among \'etale groupoids. The natural philosophy coming out from foliation theory is to consider  \'etale groupoids only up to equivalence. 

For instance if $T_0$ is an open subset of $T$ meeting every orbit of $\Cal G$, then the inclusion in $\Cal G$ of the restriction of $\Cal G$ to $T_o$ is an equivalence. If $\Cal U = \{U_i\}_{i \in I}$ is an open cover of $T$, the localization of $\Cal G$ over $\Cal U$ is the \'etale groupoid $\Cal G_\Cal U$ whose space of objects is the disjoint union $T_\Cal U$ of the $U_i$, i.e the space of pairs $(i,x)$ with $x \in U_i$. The  morphisms are the triple $(j,g,i)$ such that $s(g) \in U_i$ and $t(g) \in U_j$. The source and target of $(j,g,i)$ are $(i,s(g))$ and $(j,t(g))$ and the composition $(k,g',j)(j,g,i)$ whenever defined is equal to $(k,g'g,i)$. The natural projection $(j,g,i) \to g$ is an equivalence from $\Cal G_\Cal U$ to $\Cal G$. Note that two \'etale groupoids $\Cal G$ and $\Cal G'$ are equivalent iff there are open cover $\Cal U$ of $T$ and $\Cal U'$ of $T'$ such that $\Cal G_\Cal U$ is isomorphic to $\Cal G'_{\Cal U'}$.

\subheading{3. Developability} A connected \'etale groupoid $\Cal G$ is developable if it is equivalent to a groupoid $\Gamma \ltimes X$ given by  a discrete group $\Gamma$ acting by homeomorphisms on the space $X$. If this is so, we can always assume that $X$ is simply connected. In that case the groupoid $\Gamma \ltimes X$ is unique up to isomorphism: $\Gamma$ is isomorphic to the fundamental group of $\Cal G$ and $X$ is equivalent to its universal covering. This is  specific to \'etale groupoids.

\subheading{4. Groupoid of germs} Another feature specific to \'etale groupoids is that we can associate to an \'etale groupoid $\Cal G$ with space of objects $T$ its \'etale groupoid of germs $\Cal H$ constructed as follows. Given $g \in \Cal G$ with source $x$,  let $\tilde g: U \to \Cal G $ be a local (continuous) section of the source projection $s$ defined on an open  neighbourhood $U$ of $s (g)$ such that $s(\tilde g(x))=g$. The germs of $t\circ \tilde g: U \to T$ at the various points of $U$ form an open set in $\Cal H$. The map $g \mapsto$ germ at $x$ of $t \circ \tilde g$ is an \'etale  surjective homomorphism $\Cal G \to \Cal H$. A geometric structure on $T$, as for instance a differentiable or analytic or Riemannian structure, invariant by $\Cal H$ induces the corresponding geometric structure on $\Cal G$.

\heading{II. Morphisms and cocycles}\endheading

\subheading{1. $(\Cal G',\Cal G)$-bundles}
Let $\Cal G$ and $\Cal G'$ be \'etale groupoids with space of objects $T$ and $T'$ respectively. We recall the notion  of $(\Cal G',\Cal G)$-bundle sketched to me by Georges Skandalis after my Toulouse talk in 1982 (see \cite 5 for a general definition in case of  topological groupoids and \cite 7 for differentiable groupoids). A $(\Cal G',\Cal G)$-bundle is a topological space $E$ together with  a right action of $\Cal G$ with respect to a continuous projection $s:E \to T$ and a left action of $\Cal G$ with respect to a continuous map $t:E \to T'$. The following conditions must be satisfied.

a) $E$ is a left principal $\Cal G'$-bundle with base space $T$. This means that $s:E \to T$ is surjective and that each point of $T$ has an open neighbourhood $U$ with a continuous section $\sigma: U \to E$ (i.e. $s\circ \sigma$ is the identity of $U$) such that the map $\Cal G' \times_{T'}U \to s^{-1}(U)$ sending $(g',x)$ to $g'.\sigma(x)$ is an isomorphism. Here $\Cal G' \times_{T'} U$ is the subspace of $\Cal G' \times  U$ formed by the pairs $(g',x)$ with $s(g') = t(\sigma(x))$.

b) The right action of $\Cal G$ commutes with the left action of $\Cal G'$.

\medskip

Let $\Cal G''$ be an \'etale groupoid with space of objects $T''$ and let $E'$ be a $(\Cal G'',\Cal G')$-bundle. Then the composition $E'\circ E$ is the $(\Cal G'',\Cal G)$-bundle $E'\times _{\Cal G'} \Cal G$-bundle which is the quotient of $E' \times_{T'} E$ by the equivalence relation which identifies $(e'.g',e)$ to $(e',g'.e)$ for $g' \in \Cal G'$ with $s(e') = t(g')$ and $s(g') = t(e)$. The equivalence class of $(e',e)$ is noted $[e',e]$. The projections source and target map $[e',e]$ to $s(e)$ and $t(e')$. The actions of $g'' \in \Cal G''$ and $g \in  \Cal G$ are given by $g''.[e',e].g = [g''.e',e.g]$.

\remark{2. Example} Let $\phi:\Cal G \to \Cal G'$ be a homomorphism inducing the map $\phi_0:T \to T'$ on the spaces of objects. The associated $(\Cal G',\Cal G)$-bundle $E_\phi$ is the space $E_\phi = \Cal G' \times_{T'} T$ made up of pairs $(g',x) \in \Cal G' \times T$ such that $s(g') =\phi_0(x)$. The source (resp. target) map sends $(g',x)$ to $x$ (respectively $t(g')$). The action of $ g_1' \in \Cal G'$ with $s( g_1')= t(g')$ on $(g',x)$ is given by $g'_1.(g',x) = (g'_1g',x)$. The action of $g \in  \Cal G$ with $t(g) =x$ and $s(g) = y$ on $(g',x)$ is given by $(g',x).g = (g'\phi(g),y)$. The projection $s:E_\phi \to T$ has a global section $\sigma:T \to E_\phi$ given by $\sigma(x) = (\phi(1_x),x)$. If $\phi':\Cal G' \to \Cal G''$ is a homomorphism, then $E_{\phi'}\circ E_\phi$ is naturally isomorphic to $E_{\phi' \circ \phi}$.

 Conversely if a $(\Cal G',\Cal G)$-bundle $E$ has a global section $\sigma:X \to E$, then it is canonically isomorphic to $E_\phi$, where $\phi:\Cal G \to \Cal G'$ is the homomorphism defined by $$\sigma(t(g)).g =\phi(g).\sigma(s(g)).$$ Another section $\sigma':T \to E$ gives a homomorphism  related to $\phi$ by the natural transformation $h$ defined by $\sigma'(x)= h(x)\sigma(x)$.

\endremark

\subheading{3. Morphisms and equivalences} An isomorphism between two $(\Cal G',\Cal G)$-bundles $E$ and $E'$ is a homeomorphism from $E$ to $E'$ commuting with $s$ and $t$  and which is $\Cal G'$ and $\Cal G$-equivariant. The  isomorphism class of a $(\Cal G',\Cal G)$-bundle $E$ is noted $[E]$ and is called a morphism from $\Cal G$ to $\Cal G'$. The set of morphisms from $\Cal G$ to $\Cal G'$  is noted $\Cal Mor(\Cal G',\Cal G)$. 

A $(\Cal G',\Cal G')$-bundle is invertible if it is also a right principal  $\Cal G$-bundle over $T'$ with respect to the projection $t$. It is also called a $(\Cal G',\Cal G)$-bibundle. The inverse $E^{-1}$ of $E$ is the $(\Cal G,\Cal G')$-bundle obtained from $E$  by exchanging the role of $s$ and $t$ and the actions by their opposite. Note that $E \times_{\Cal G} E^{-1}$ is isomorphic to $\Cal G'$ and $E^{-1}\times_{\Cal G'} E$ is isomorphic to $\Cal G$. The isomorphism class $[E]$ of an invertible bundle is called an equivalence from $\Cal G$ to $\Cal G'$. The bundle associated to a homomorphism $\phi$ represents an equivalence iff $\phi$ is an equivalence.

The morphisms from \'etale groupoids to \'etale groupoids form a category. Note that if $T$ and $T'$ have invariant geometric structures, one can also consider morphisms and equivalences compatible with those structures.

\subheading{4. Description in terms of cocycles} Given a $(\Cal G',\Cal G)$-bundle $E$, let $\Cal U =\{U_i\}_{i \in I}$ be an open cover of $T$ such that there are continuous sections $\sigma_i;U_i \to E$ of $s$ for each $i \in I$. To this  we can associate the homomorphism $\phi: \Cal G_\Cal U\to \Cal G'$ mapping $(j,g,i)$, with $s(g)=x$ and $ t(g)=y$,  to the unique element $\phi(j,g,i)$ such that
$$ \sigma_j(y)g = \phi(j,g,i)\sigma_i(x) .$$

Conversely, a homomorphism $\phi: \Cal G_\Cal U \to \Cal G'$ gives a $(\Cal G',\Cal G)$-bundle $E_{\phi,\Cal U}$ as follows. Let $\Cal G' \times_{T'}U_i$ be the subspace of $\Cal G' \times U_i$ consisting of pairs $(g',x)$ such that $\phi_0(x,i) = s(g')$. In the disjoint union $\coprod_{i \in I} \Cal G' \times_{T'}{U_i}$ consider the equivalence relation identifying, for $x \in U_i \cap U_j$, the element $(j,g',x) \in \{j\} \times\Cal G' \times_{T'}{U_j} $ to the element $(i,g'\phi(j,1_x,i),x) \in \{i\} \times\Cal G' \times_{T'}{U_i} $. The equivalence class of $(i,g',x)$ is noted $[i,g',x]$ and the quotient of 
$\coprod \Cal G' \times_{T'}{U_i}$ by this equivalence relation is noted $E_{\phi,\Cal U}$. The projections $s$ and $t$   map $[i,g',x]$ to $x$ and $t(g')$ respectively. For  $g'_1 \in \Cal G' $ with $s(g'_1) = t(g')$ we define $g'_1.[i,g',x] = [i,(g'_1g').x]$. For $g \in \Cal G$ wiht source $y \in U_j$ and target $x \in U_i$, we define $[i,g',x].g = [j,g'\phi(i,g,j),y]$. This defines on $E_{\phi,\Cal U}$ the structure of a $(\Cal G',\Cal G)$-bundle, called the $(\Cal G',\Cal G)$-bundle constructed from $\phi$.

Let $\Cal U' =\{U'_{i'}\}_{i' \in I'}$ be another open cover of $T$ and $\phi': \Cal G_{\Cal U'} \to \Cal G'$ be a homomorphism. We assume that the sets of indices $I$ and $I'$ are disjoint.  Let $\Cal U \coprod \Cal V$ be the union of these open covers. We can identify $\Cal G_\Cal U$ and $\Cal G_{\Cal U'}$ to disjoint subgroupoids of $\Cal G_{\Cal U \coprod \Cal U'}$. The $(\Cal G',\Cal G)$-bundles $E_{\phi,\Cal U}$ and $E_{\phi',\Cal U'}$ associated to $\phi$ and $\phi'$ are isomorphic iff there is a homomorphism $\psi: \Cal G_{\Cal U \coprod \Cal U'} \to \Cal G'$ extending $\phi$ and $\phi'$. Such an extension gives a well defined isomorphism between these bundles.
\medskip

\heading{III. Pointed morphisms}\endheading

\definition{1. Definition} A pointed $(\Cal G',\Cal G)$-bundle over a point $* \in T$ (see \cite {10}) is a pair $(E,e_0)$, where $E$ is a 
$(\Cal G',\Cal G)$-bundle and $e_0 \in E$ a point such that $s(e_0) = *$. A pointed morphism from $\Cal G$ to $\Cal G'$ over $*$ is the isomorphism class $[E,e_0]$ of a pointed bundle $(E,e_0)$ over $*$. The set of pointed morphisms over $*$ is noted $\Cal Mor(\Cal G',\Cal G,*)$. In the case where $\Cal G$ is the circle (considered as a trivial groupoid), a pointed morphism from $\Cal G$ to $\Cal G'$ was called in \cite 4  a based $\Cal G'$-loop. \enddefinition

Note that the bundle $E_\phi =\Cal G' \times_{T'} T$ associated to a homomorphism $\phi: \Cal G \to \Cal G'$ inherits naturally a base point over $* \in T$, namely the point $(\phi(1_*),*)$. Also the bundle $E_{\phi,\Cal U}$ associated to a homomorphism $\phi: \Cal G_{\Cal U} \to \Cal G'$ inherits a base point once we have chosen $U \in \Cal U$ containing $*$, namely the class of the point $(\phi(1_*),*) \in \Cal G' \times_{T'} U$.

$\Cal G'$ acts on the left on $\Cal Mor(\Cal G',\Cal G,*)$ with respect to the map $t:\Cal Mor(\Cal G',\Cal G,*) \to T'$
sending $[E,e_0]$ to $t(e_0)$. The action of $g'\in \Cal G'$ on $[E,e_0]$ is defined by $g'.[E,e_0] = [E,g'.e_0]$. The associated groupoid  with set of objects $\Cal Mor(\Cal G',\Cal G,*)$ will be noted $\Cal G' \times_{T'}\Cal M(\Cal G',\Cal G,*)$. We shall introduce below, under suitable conditions, a topology  which makes it an \'etale  groupoid. For that the following lemma will be crucial. Note that the set of orbits $\Cal G' \backslash \Cal Mor(\Cal G',\Cal G,*)$ is isomorphic to $\Cal Mor(\Cal G',\Cal G)$.

\proclaim{2. Lemma}Suppose that $\Cal G$ is connected and that $\Cal G'$ is Hausdorff. Then any automorphism $h$ of $E$ preserving a base point $e_o \in E$ is the identity. 
\endproclaim
\demo{Proof}  Let $F \subset E$ be the set of points $e \in E$ such that $h(e) =e$. The set $F$ is open because it preserves the projection $s$ which is \'etale. The set $F$ is closed because $\Cal G'$ is Hausdorff. It is non empty and  invariant under $\Cal G$ and $\Cal G'$. Hence it is equal to $E$ because $(\Cal G' \times \Cal G)\backslash E = \Cal G\backslash T$ is assumed to be connected. $ \square$ \enddemo

In terms of cocycles, with the notations at the end of II.4, the lemma is formulated as follows. Assume that we have chosen $i \in I$ and $i' \in I'$ such that $* \in U_i \cap U'_{i'}$. Let $(E_{\phi,\Cal U},e_0)$ and $(E_{\phi',\Cal U'},e'_0)$ be the pointed bundles associated as above to $\phi$ and $\phi'$. Then these two pointed bundles are isomorphic iff there is a homomorphism $\psi: \Cal G_{\Cal U \coprod \Cal U'} \to \Cal G'$ extending $\phi$ and $\phi'$ such that $\psi(i,1_*,i')=\phi(i,1_*,i) = \phi'(i',1_*,i')$. Under the hypothesis of the lemme, the assertion is that such a $\psi$ is unique.

\proclaim{3. Proposition} For $\Cal G$ connected , the equivalence class of $\Cal G \times _{T'} \Cal Mor(\Cal G',\Cal G,*)$ does not change if we replace $\Cal G$ and $\Cal G'$  by equivalent \'etale groupoids, or if we change the base point $* \in T$.\endproclaim 
\demo{Proof} For instance, assume that $\Cal E'$ is a bibundle giving an  equivalence from $\Cal G'$ to an \'etale groupoid $\overline {\Cal G}'$ with space of objects  $\overline T'$. 
Then an equivalence from   $ \Cal G' \times_{T'} \Cal Mor( \Cal  G',\Cal G,*)$ to $\overline {\Cal G}' \times_{T'} \Cal Mor(\overline {\Cal  G}',\Cal G,*)$ will be represented by the bibundle $\Cal E' \times_{T'} \Cal Mor(\Cal G',\Cal G,*)$. 
The projection $s$ maps  $(e',[E,e_0])$ to $[E,e_0]$ and the projection $t$ maps it to $(\Cal E' \times_{\Cal G'}E, [e',e_0])$. An element of $\Cal G' \times_{T'} \Cal Mor(\Cal  G',\Cal G,*)$ with target $[E.e_0]$ is of the form $(g',[E,{g'}^{-1}e_0])$, where $g'\in \Cal G'$ with $t(g') =t(e_0)$. Its right action on $(e',[E,e_0])$ is equal to $(e'g',[E,{g'}^{-1}e_0])$. An element of $\overline {\Cal G}' \times_{T'} \Cal Mor(\overline{\Cal  G}',\Cal G,*)$ with source $[\Cal E'\times_{\Cal G'}E, [e',e_0])$ is of the form $(\overline g',[\Cal E'\times_{\Cal G'}E, [e',e_0])$, where $s(\overline g') = t(e')$. Its left  action on $(e',[E,e_0])$ is equal to $(\overline g'e',[E,e_0])$. With those definitions, it is easy to check that $\Cal E'\times_{T'}\Cal Mor(\Cal G',\Cal G,*)$ is a bibundle giving an equivalence from $\Cal G' \times_{T'}\Cal Mor(\Cal G',\Cal G,*)$ to $\overline{\Cal G'}\times_{\overline T'}\Cal Mor(\overline  {\Cal G}',\Cal G,*)$.

The other claims of the proposition are similarly easy to prove. For instance, if $*' \in T$ is another base point, the elements of the bibundle giving the equivalence from $\Cal G'\times_{T'}\Cal Mor(\Cal G',\Cal G,*)$ to $\Cal G' \times_{T'}\Cal Mor(\Cal G',\Cal G,*')$ are the isomorphism classes of bipointed $(\Cal G',\Cal G)$-bundles over $*$ and $*'$. $\square$ \enddemo

As an example, we consider below the developable case.

\proclaim{4. Proposition} Assume that $\Cal G = \Gamma \ltimes T$ and $\Cal G' = \Gamma' \ltimes T'$, with $T$ simply connected. 

Then 
$\Cal Mor(\Cal G',\Cal G,*) $ is the set of pairs $(f,\psi)$, where $\psi:\Gamma \to \Gamma'$ is a group homomorphism and $f:T \to T'$ is a $\psi$-equivariant continuous map.
 
The groupoid $\Cal G' \times_{T'}\Cal Mor(\Cal G',\Cal G,*)$ is isomorphic to $\Gamma' \ltimes \Cal Mor(\Cal G',\Cal G,*)$, where the action of $\gamma'$ on the pair $(f,\psi)$ is given by
$$\gamma'.(f,\psi) = (t_{\gamma'}\circ f, Ad(\gamma')\circ \psi),$$
and $t_{\gamma'}$ is the translation of $T'$ by $\gamma'$. \endproclaim

\demo{Proof} Let $(E,e_o)$ be a pointed bundle representing an element of $\Cal Mor(\Cal G',\Cal G,*)$. With respect to the projection $s:E \to T$, the bundle  $E$ is a $\Gamma'$-principal covering. As $T$ is simply connected, there is a unique lifting  $\sigma:T \to E$ of the projection $s$ such that $\sigma(*) = e_0$. Then to $\sigma$ is associated as in II.2 a continuous homomorphism $\phi: \Gamma \ltimes T \to \Gamma' \ltimes T'$ such that $E$ is canonically isomorphic to $E_\phi$.  Let  $f:T \to T'$ be the map $t \circ \sigma$  induced by $\phi$ on the space of objects. As $T$ is connected and $\Gamma'$ is discrete,   $\phi(\gamma,x)$ is of the form $(\psi(\gamma), f(x))$. Composing with the target projection, we get $f(\gamma.x)=\psi(\gamma).f(x)$. Conversely such a pair $(f,\psi)$ defines a homomorphism $(\gamma,x) \mapsto (\psi(\gamma),f(x))$ from $\Cal G$ to $\Cal G'$.

The choice of another base point over $*$ leads to a homomorphism $\phi'$ related to $\phi$ by a natural transformation $h:T\to\Cal G'$ of the form $h(x) =(\gamma', f(x))$. So $\phi'(\gamma,x)= (\gamma'\psi(\gamma){\gamma'}^{-1}, \gamma'.f(x))$. $\square$ \enddemo

\subheading{5. The exponential morphism} We want to define a morphism from the groupoid product of $\Cal G' \times_{T'} \Cal Mor(\Cal G', \Cal G,*)$ with $\Cal G$ to the groupoid $\Cal G'$, which is the usual exponential map ${T'}^T \times T \to T'$ when the groupoids $\Cal G$ and $\Cal G'$ are the trivial groupoids $T$ and $T'$. We first describe the principle of the construction on the set theoretical level, under the hypothesis of the lemma.  In IV we shall consider topologies.

For each $z \in Mor(\Cal G', \Cal G,*)$, choose a pointed  $(\Cal G',\Cal G)$-bundle $(E^z,e^z_0)$ whose isomorphism class  $[E^z,e^z_0]$ is $z$. Let $EXP$ be the  disjoint union
$\coprod_{z \in Mor(\Cal G', \Cal G,*)} E^z$. We have projections $s:EXP \to Mor(\Cal G', \Cal G,*)\times T$ and $t:EXP \to T'$ mapping $(z,e \in E^z)$ respectively to $(z,s(e))$ and $t(e)$. For $(g',z') \in \Cal G' \times_{T'} Mor(\Cal G', \Cal G,*)$ with target $z$ we have an  isomorphism $m_{(g'.z')}$ from the pointed bundle $(E^{z'},g'.e^{z'}_0)$ to the pointed bundle $(E^z,e^z_0)$. The lemma implies that this isomorphism is unique. Commuting left action of $\Cal G'$ and right action of $\Cal G' \times_{T'} Mor(\Cal G', \Cal G,*) \times \Cal G$ on $EXP$ are defined as follows. For $(z,e) \in (z,E^z) \subset EXP$, $g' \in \Cal G'$ with source $t(e)$,$(g'',z') \in \Cal G' \times_{T'} Mor(\Cal G', \Cal G,*)$ with target $z$ and $g \in \Cal G$ with target $s(e)$ we define
$$g'.(z,e).((g'',z'),g) = (z',m_{(g'',z')}^{-1}(g'.e.g)) \in (z',E^{z'}) \subset EXP.$$
Note that $m_{(g'',z')}^{-1}(g'.e.g) =g'.m_{(g'',z')}^{-1}(e).g$. These actions define on  $EXP$ the structure of a $(\Cal G',\ (\Cal G' \times_{T'} Mor(\Cal G', \Cal G,*)) \times \Cal G)$-bundle defining the desired  exponential morphism.

\remark{Example} Assume that $\Cal G = \Gamma \ltimes T$ and $\Cal G' = \Gamma'\ltimes T'$ with $T$ simply connected as in III.4. Then $EXP$ is the disjoint union of the $(\Cal G',\Cal G)$-bundles $E^{(f,\psi)}$ indexed by the set of pairs $(f ,\psi) \in \Cal M(\Cal G',\Cal G.*)$. The bundle $E^{(f,\psi)}$ is the product $\Gamma' \times T$ with projections $s$ and $t$ mapping $(\gamma',x)$ to $x$ and $\gamma'.f(x)$ respectively. The right action of $(\gamma, \gamma^{-1}.x)\in \Gamma\ltimes T$ on $(\gamma',x)$ is equal to $(\gamma'\psi(\gamma),\gamma^{-1}.x)$ and the left action of $(\gamma'',\gamma'.f(x)) \in \Gamma'\ltimes T'$ on $(\gamma',x)$ is equal to $ (\gamma''\gamma',x)$.
In that case the bundle $EXP$ has a canonical section and therefore the exponential morphism can be described directly by the homomorphism from $(\Cal G' \times_{T'}\Cal Mor(\Cal G',\Cal G,*)) \times \Cal G= (\Gamma' \ltimes\Cal Mor(\Cal G',\Cal G,*)) \times (\Gamma \ltimes T)$  to $\Cal G=\Gamma'\ltimes T'$ mapping $(\gamma'',(f,\psi),\gamma,x)$ to $(\gamma''\psi(\gamma),f(x))$.

\heading  IV.  A topology on $\Cal Mor(\Cal G',\Cal G,*) $\endheading
To get some feeling we begin with the  particular case considered in III.4
.
\subheading{1. The developable case} We assume that 

1) $ \Cal G' = \Gamma' \ltimes T'$, where $T'$ is differentiable manifold with a $\Gamma'$-invariant Riemannian metric.

2) $\Cal G= \Gamma\ltimes T$, where $T$ is $1$-connected and there is a compact subset of $T$ meeting every orbit.

We want to define an open neighbourhood of an element $(f,\psi) \in \Cal Mor(\Cal G',\Cal G,*)$ (see III.4). Let $\tau T'$ be the tangent bundle of $T'$ and let
$f^*(\tau T') = \{(x,v), x\in T, v \in \tau_{f(x)}T'\}$ be its pull back by $f$. The group $\Gamma$ acts naturally on it, namely $\gamma.(x,v)= (\gamma.x, \psi(\gamma).v)$, where $\psi(\gamma).v$ denote the image of $v$ by the differntial of the isometry defined by $\psi(\gamma)$. Let $V_{(f,\psi)}$ be the Banach space of $\Gamma$-invariant  continuous sections of $f^*(\tau T')$ with the sup-norm. Let $\epsilon$ be a small enough positive number such that the exponential map $exp_{f(x)}:\tau_{f(x)}T' \to T'$ restricted to the open ball of radius $\epsilon$ is a diffeomorphism on a convex geodesic ball of $T'$, and this for all $x \in T$ (this is possible by condition 2)). Let $V^\epsilon_{(f,\psi)}$ be the open ball of radius $\epsilon$ in $V_{(f,\psi)}$.

We define a structure of Banach manifold on $\Cal Mor(\Cal G',\Cal G,*)$. A chart at $(f,\psi)$ is the map $V^\epsilon_{(f,\psi)} \to \Cal Mor(\Cal G',\Cal G,*)$ sending a section $\nu:x \mapsto (x,v(x))$ to $(f^\nu,\psi)$,  where $f^\nu(x) =exp_{f(x)}v(x)$. One verifies as usual  (see Eells \cite 3) that the changes of charts are differentiable. 

With this topology, the exponential homomorphism as described at the end of III is continuous.

\subheading{2. A more general case} We make the following hypothesis which is satisfied for $\Cal G'$ the transverse  holonomy or monodromy  groupoid of a Riemannian foliation on a complete  
Riemannian manifold (see \cite 9 and the appendix by E. Salem), in particular for Riemannian orbifolds.

1) $\Cal G'$ is a Hausdorff and complete Riemannian \'etale groupoid. This means that $T'$ has a $\Cal G$-invariant  Riemannian metric, that $\Cal G$ is Hausdorff and that the following condition is satisfied: for each $g' \in \Cal G'$ with source $x'$ and target $y'$ which are centers of convex geodesic balls $B(x',\epsilon)$ and $B(y',\epsilon)$, there is a unique continuous map $\tilde {g'}:B(x',\epsilon) \to \Cal G'$ such that, for each $z \in B(x',\epsilon)$, the source of $\tilde {g'}(z)$ is $z$ and $g'= \tilde{g'}(x')$. In particular $t\circ\tilde{ g'}\circ exp_{x'} = exp_{y'}\circ Dg'$, where $Dg'$ is the differential of $t\circ \tilde {g'}$ at $x'$ (it depends only on $g'$).

 It implies the following. If $g'' \in \Cal G'$ with source $y'$ and target  the center of a convex geodesic ball of radius $\epsilon$, then 
$$\widetilde {g''g'}(z) = \tilde {g''}(t(\tilde {g'}(z)))\tilde{ g'}(z)  \ \ \ \forall z \in B(x',  \epsilon).  \tag *$$

2) $\Cal G$ is connected and there is a relatively compact open subset $T_0 \subseteq T$ meeting every orbit of $\Cal G$.

\proclaim{3. Theorem} Under the above hypothesis, there is a natural structure of Banach manifold on $\Cal Mor(\Cal G',\Cal G,*)$ so that $\Cal G' \times_{T'}\Cal M(\Cal G',\Cal G,*)$ becomes a differentiable Banach \'etale groupoid. \endproclaim

\demo{Proof} The proof will be in three steps. 

i) Given a pointed $(\Cal G',\Cal G)$-bundle $(E,e_0)$ over $* \in T$, we describe a chart at $[E,e_0]$. Let $\Cal G_0$ be the restriction of $\Cal G$ to $T_0$.  The restriction $E_0$ of $E$ above $T_0$ is a pointed $(\Cal G',\Cal G_0)$-bundle and $(E,e_0)$ is naturally isomorphic to $E_0\times_{\Cal G_0} \Cal E$, where $\Cal E =\{g \in \Cal G: t(g) \in T_0\}$, with the base point $[e_0,1_*]$.

Choose a finite open cover $\Cal U$ of $T_0$ consisting of $U_0,U_1.\dots,U_k$ such that $* \in U_0$ and such that there are continuous sections $\sigma_i:U_i \to E$ with $\sigma_0(*)=e_0$. Let $\Cal G_\Cal U$ be the localization of $\Cal G_0$ over  $\Cal U$. As in II.4, we get a homomorphism $\phi:\Cal G_{\Cal U} \to \Cal G'$ whose restriction to the space of objects is a continuous map $\phi_0: T_\Cal U \to T'$. The groupoid $\Cal G_\Cal U$ acts on $\phi_0^*(\tau T')$: the action of $g \in \Cal G_\Cal U$ with source $x$ and target $y$  on $(x,v(x))$, where $ v(x) \in \tau_{\phi_0(x)} T'$, is defined by $g.(x,v(x)) = (y, D\phi(g).v(x))$. Let $V_\phi$ be the Banach space of continuous $\Cal G_\Cal U$-invariant sections of $\phi_0^*( T')$. Thanks to hypothesis 2), we  can find a positive number $\epsilon$ such that, for all  $x \in T_\Cal U$, the exponential map $exp_{\phi_0(x)}$ gives a diffeomorphism from the $\epsilon$-ball in $\tau_{\phi_0(x)}T'$ to a convex geodesic ball in $T'$. Let $V^\epsilon_\phi$ be the open subset of $V_\phi$ formed by the $\Cal G_\Cal U$-invariant sections of $\phi_0^*( T')$ of norm $<\epsilon$.

Given $g \in \Cal G_\Cal U$ with source $x$, let $\widetilde {\phi(g)}$ be the extension of $\phi(g)$ to the geodesic ball of radius $\epsilon$ centered at the source $\phi_0(x)$ of $\phi(g)$ as defined in 1). Given a section $\nu:x \mapsto (x,v(x))$ in $V^\epsilon_\phi$, let $\phi^\nu:\Cal G_\Cal U \to\Cal G'$ be defined, for $g$ with source $x$, by 
$$\phi^\nu(g) = \widetilde {\phi(g)}(exp_{\phi_0(x)}v(x)).$$
From $\phi^\nu$, we reconstruct as in II.4 a $(\Cal G',\Cal G_0)$-bundle $E^\nu_0 = E_{\phi^\nu,\Cal U}$ with base point $e^\nu_0$ (see III.1) and then, as above, a $(\Cal G',\Cal G)$-bundle $E^\nu$ with a base point still noted $e^\nu_0$.

The map $\nu \mapsto [E^\nu,e^\nu_0]$ from  
$V^\epsilon_\phi$ to $\Cal Mor(\Cal G',\Cal G,*)$ will be a chart $h_\phi$ for the structure of Banach manifold. This map is injective because for $\nu, \nu' \in V^\epsilon_\phi$, the set of points $x$ in the space of objects of $\Cal G_\Cal U$ such that $v(x) =v'(x)$ is a closed $\Cal G_\Cal U$-invariant subset which is also open. If $h_\phi (\nu) =h_\phi(\nu')$ it is non empty because it contains the base point, hence it is the whole of the space of objects  because $\Cal G_\Cal U$ is connected.

ii) Let  $\phi': \Cal G_{\Cal U'} \to \Cal G'$ be another  pointed homomorphism, where  $\Cal U'$ is another open cover $\{U'_0,\dots, U'_{k'}\}$  of $T_0$ such that $* \in U'_0$ , and let $h_{\phi'}: V^{\epsilon'}_{\phi'} \to \Cal Mor(\Cal G',\Cal G,*)$ be a corresponding chart. 
One has to check that the change of charts $h_{\phi'}^{-1}{h_\phi}^{-1}$
is differentiable. We first consider some particular cases.
 
a) Assume that $\phi :\Cal G_\Cal U \to \Cal G'$ and $\phi':\Cal G_{\Cal U'}\to \Cal G'$ are pointed equivalent (see III,2). Then there is a naturel linear isomorphism from $V_\phi$ to $V_{\phi'}$ and the change of charts $h_{\phi'}^{-1} h_\phi$ is the restriction of this isomorphism to the balls of radius the minimum of $\epsilon$ and $\epsilon'$

b) If $\phi' = \phi^\nu  :  \Cal G_\Cal U \to \Cal G'$, where $\nu \in V^\epsilon_\phi$, the usual argument shows that the change of charts  $h_{\phi'}^{-1} h_\phi$ is differentiable at $\nu$.

In the general case, if $h_\phi(\nu) = h_{\phi'} ({\nu'})$, we express the change of charts $h_{\phi'}^{-1} h_\phi$
around $\nu$ as the composition $(h^{-1}_{\phi'}h_{{\phi'}^{\nu'}})(h^{-1}_{{\phi'}^{\nu'}}h_{\phi^\nu})(h_{\phi^\nu}^{-1} h_\phi)$ and we apply a) and b).

The map $\Cal M(\Cal G',\Cal G,*) \to T'$ sending $[E,e_0]$ to $t(e_0)$ is continuous because, for each chart $h_\phi$, the map $\nu \mapsto t(e^\nu_0)
=exp_{\phi_0(*)}(v(*))$ is continuous.

iii) Let us check now that  the action of $\Cal G'$ on $\Cal Mor(\Cal G',\Cal G,*)$ is differentiable. Let $(E,e_0)$ be a pointed bundle and let $g' \in \Cal G'$ be such that $
s(g') = t(e_0)$. Choose an open neighbourhood $U_0$ of $*$ in $T_0$ such that there are continuous sections $\sigma_0, \sigma'_0: U_0 \to E$ such that $\sigma_0(*) =e_0$ and $\sigma'_0(*)=g'.*$. We can find a finite open cover $\Cal U =\{U_0, \dots, U_k\}$ of $T_0$ such that there are continuous sections $\sigma_i:U_i \to E$ for $i \neq 0$. Let $\phi, \phi': \Cal G_\Cal U \to \Cal G'$ be the homomorphisms associated to the choice of the sections $\sigma_i$ and $\sigma'_i$, where $\sigma'_i = \sigma_i$ for $i \neq 0$. The associated pointed bundles are naturally isomorphic to $(E,e_0)$ and $(E,g'.e_0)$. For $\epsilon > 0$ small enough, the charts $h_\phi: V^\epsilon_\phi \to \Cal Mor(\Cal G',\Cal G,*)$ and $h_{\phi'}: V^\epsilon_{\phi'} : \Cal Mor(\Cal G',\Cal G,*)$ are defined. For each $x \in U_0$, let $g'(x)$ be the unique element of $\Cal G'$ such that $\sigma'(x) = g'(x).\sigma(x)$. 
The extension $\widetilde{g'(x)}$ of $g'(x)$ is defined on the geodesic ball of radius $\epsilon$ and center $s(g'(x)$. Let $f:V_\phi \to V_{\phi'}$ be the linear isometry mapping the $\Cal G_\Cal U$-invariant section $\nu: (i,x) \in T_\Cal U \mapsto ((i,x), v(i,x))$ in $V_\phi$ to the $\Cal G_{\Cal U'}$- invariant section $f(\nu) \in V_{\phi'}: (i,x) \mapsto ((i,x),f(v)(i,x))$, where $f(v)(i,x)$ is equal to $v(i,x)$ for $i \neq 0$ and $f(v)(0,x) = Dg'(x). v(0,x)$. Let $\sigma$ be the section of $s: \Cal G' \times_{T'}\Cal Mor(\Cal G',\Cal G.*) \to \Cal Mor(\Cal G',\Cal G,*)$ above $h_\phi(V^\epsilon_\phi)$ given by $h_\phi(\nu) \mapsto (\widetilde { g'}(t(e_0^\nu),h_{ \phi}(\nu))$. The differentiable structure on $\Cal M(\Cal G',\Cal G,*)$ is invariant by $\Cal G'$ because   $h_{\phi'}^{-1}\circ t \circ \sigma \circ h_\phi = f$ is differentiable.
$\square$\enddemo

\subheading{4. The exponential morphism} The canonical morphism $EXP$ from $(\Cal G' \times_{T'}\Cal Mor(\Cal G',\Cal G,*)) \times \Cal G$ to $\Cal G'$ is represented in the above framework by a tautological bundle over $\Cal Mor(\Cal G',\Cal G,*) \times T$ described as follows. As $\Cal Mor(\Cal G',\Cal G,*)$ is canonically isomorphic to $\Cal Mor(\Cal G',\Cal G_0,*)$, we can replace $\Cal G$ by it s rectriction $\Cal G_0$ to $T_0$ (see the beginning of the proog of 3). To each chart $h_\phi$ we associate the family of pointed $(\Cal G',\Cal G)$-bundles
$(E_{\phi^\nu,\Cal U}, e_0^\nu)$ parametrized by the elements $\nu \in V^\epsilon_\phi$. We note $EXP_\phi$ the disjoint union $\coprod_{\nu \in V^\epsilon_\phi} E_{\phi^\nu,\Cal U}$ topologized as the quotient of the subspace of $V^\epsilon_\phi \times\Cal G' \times T_\Cal U$ made up of the triples $(\nu,g',(i,x))$ such that $s(g') = \phi_0^\nu(i,x)$. We have a continuous projection $s_\phi: EXP_\phi \to V^\epsilon_\phi \times T_0$  sending $(\nu,e)$ to $(\nu,s(e))$. 

Let $h_{\phi'}: V^{\epsilon'}_{\phi'} \to \Cal Mor(\Cal G',\Cal G,*)$ be another chart. The change of charts $h_{\phi'}^{-1}h_\phi: W \subset V^\epsilon_\phi  \to W' \subset V^{\epsilon'}_{\phi'}$ 
lifts functorially to a family of pointed isomorphisms $f_{\phi', \phi}: s_{\phi}^{-1}(W \times T_0) \to s_{\phi'}^{-1} (W' \times T_0)$. The bundle over $\Cal Mor(\Cal G',\Cal G,*) \times T_0$ representing $EXP$ is the quotient of the disjoint union of 
the $EXP_\phi$ under the equivalence relation which identifies points corresponding to each other by the isomorphisms $f_{\phi' \phi}$. This bundle is also denoted $EXP$ and its restriction above the range of the chart $\phi$ is canonically isomorphic to $EXP_\phi$. The map $f_\phi: E_\phi \to EXP$ associating to a point its equivalence class is injective and can be considered as a chart. The change of charts $f_{\phi'}^{-1} \circ f_\phi$ is the map $f_{\phi'\phi}$. 

The left action of $\Cal G'$ and the commuting right action of $(\Cal G' \times \Cal Mor(\Cal G',\Cal G_0,*)) \times \Cal G_0$ on $EXP$ is described in the charts like in III.4.

\medskip

 Let $\Cal H$ be an \'etale groupoid with space of objects a topological space $S$. Let $\Psi$ be a morphism from $\Cal H$ to $\Cal G' \times_{T'} \Cal M(\Cal G',\Cal G,*)$.  Let $\Psi \times id_\Cal G$ be the morphism from $\Cal H \times \Cal G$ to $(\Cal G' \times_{T'}\Cal M(\Cal G',\Cal G,*)) \times \Cal G$ which is the cartesian product of $\Psi$ with the identity of $\Cal G$. 
 
 \proclaim{5. Theorem} The map associating to $\Psi$ the morphism $\overline \Psi = EXP \circ (\Psi \times id_\Cal G)$ 
 induces  a bijection between the set of morphisms from $\Cal H$ to $\Cal G' \times_{T'} \Cal Mor(\Cal G',\Cal G,.*)$ and the set of morphisms from $\Cal H \times \Cal G$ to $\Cal G'$. \endproclaim 
 
 \demo{Proof} We just indicate how starting from a bundle $\overline P$ representing a morphism $\overline \Psi$ from $\Cal H \times \Cal G$ to $\Cal G'$ we can construct a morphism $\Psi$ from $\Cal H$ to $\Cal G' \times_{T'} \Cal M(\Cal G',\Cal G,*)$. For each $v \in S$, the pull back of $\overline P$ by the inclusion $T \to \{v\}\times T $ is a $(\Cal G',\Cal G)$-bundle over $T$ noted $\overline P_v$. An element $h \in \Cal H$ with source $v$ and target $w$ gives an isomorphism  $\overline P_h:\overline P_w \to \overline P_v$ induced by the  right action of $(h,1_T)$ on $\overline P$.
 
 Choose an open cover $\Cal V =\{V_i\}_{i\in I}$ of $S$ such that there exist  continuous sections $\sigma_i: V_i \times \{*\} \to \overline P$ of the projection $\overline P \to S \times T$. For $v \in V_i$, this induces a base point above $*$ on $\overline P_v$ noted $\sigma_i(v)$. We note $f(v,i) \in \Cal M(\Cal G',\Cal G,*)$ the isomorphism class of the pointed bundle $(P_v,\sigma_i(v))$.
 
   Let $h$ be an element of $\Cal H$ with source $v \in V_i$ and target $w \in V_j$. Let $f(j,h,i)$ be the element of $\Cal G'$ such that 
 $$\sigma_j(w,*).(h,1_*) = f(j,h,i).\sigma_i(v,*).$$
 Then the map  $\psi:\Cal H_\Cal V \to \Cal G'\times_{T'} \Cal M(\Cal G'\Cal G,*)$ sending $(j,h,i) \in \Cal G_\Cal V$ to $(f(j,h,i), f(v,i))$ is a homomorphism representing the morphism $\Psi$. $\square$

 \enddemo
 
 \remark{Remark} A pointed morphism from $\Cal H$ to $\Cal G' \times_{T'} \Cal Mor(\Cal G',\Cal G,*)$       above a point $v_0 \in S$ is determined by the choice of a pointed $(\Cal G',\Cal G)$-bundle $(E,e_0)$ over $*$. It corresponds to a morphism from $\Cal H \times \Cal G$ to $\Cal G'$ represented by a $(\Cal G', \Cal H \times \Cal G)$-bundle $\overline P$ over $S \times T$ together with an isomorphism from $(E,e_0)$ to $\overline P_{v_0}$. (See \cite 4, 2.2.3 for a general notion of relative morphisms.)

\subheading{6. Remarks} If $\Cal G$ is replaced by an equivalent groupoid $\overline {\Cal G}$ and $\Cal G'$ by an equivalent groupoid $\overline {\Cal G'}$ with an invariant Riemannian metric, then (see III.2)  the groupoids $\Cal G' \times_{T'} \Cal Mor (\Cal G',\Cal G,*)$ and $\overline {\Cal G'} \times_{T'} \Cal Mor (\overline {\Cal G'},\overline {\Cal G},\overline*)$ are differentiably equivalent.

According to the needs, one can replace in the construction
of section 2 the set of continuous pointed morphisms from $\Cal G$ to $\Cal G'$ by a subset with a suitable topology, as it is usual when dealing with functions spaces (see for instance \cite 4).

\heading V. Selfequivalences and extensions \endheading

\subheading{1. The case of developable groupoids } Let $\Cal G = \Gamma \ltimes T$, where $T$ is simply connected and let $*$ be a point of $T$. According to III.4. the set $(\Cal Self(\Cal G),*)$ of pointed selfequivalences of $\Cal G$ is isomorphic to the set $\Cal S$ of  pairs $(f,\psi)$, where $\psi$ is an automorphism of $\Gamma$ and $f$ is a homeomorphism of $T$ which is $\psi$-equivariant. This set is a group, the composition being defined  by $(f,\psi)(f',\psi') = (f\circ f',\psi \circ \psi')$. The groupoid of pointed selfequivalences of $\Cal G$ is isomorphic    to  the action groupoid $\Gamma \ltimes \Cal S$, where $\gamma \in \Gamma$ acts on a pair by the rule $$\gamma.(f,\psi) = (t_\gamma \circ f,Ad(\gamma) \circ \psi).$$

In fact we have a {\bf crossed mdule}:
$$\mu: \Gamma \to \Cal S,$$
where $\mu(\gamma) = (t_\gamma,Ad(\gamma))$ and the pair $(f,\psi) \in \Cal S$ acts on $\Gamma$ through the automorphism $\psi$. 

If $T$ has a $\Gamma$-invariant Riemannian metric  and if there is a compact subset meeting every orbit, then $\Cal S$ has a topology so that $\mu$ is a {\bf topological crossed module}.

To illustrate the connection with extensions, we consider two particular cases.

\subheading{2. The case of a discrete group $\Gamma$}Here $\Gamma$ is considered as a discrete groupoid with one object $*$.  The group of pointed selfequivalences is the group $Aut(\Gamma)$ of automorphisms of $\Gamma$ and the crossed module is 
$$\mu: \Gamma \to Aut(\Gamma).$$
where $\mu$ sends $\gamma$ to $Ad(\gamma)$.

As it is well known since the works of J. H. C. Whitehead and S. Maclane in the forties, this crossed module plays the central role in the problem of extensions of groups, and more generally of extensions of groupoids by the discrete group $\Gamma$. 

Given an \'etale groupoid $\Cal G$ with space of objects $T$, an extension of $\Cal G$ by $\Gamma$ is given by an open cover $\Cal U$ of $T$, an etale groupoid $\widetilde{\Cal G}$ with space of objects $T_\Cal U$ (notation of I.2)  and a homomorphism $\phi:\widetilde{\Cal G} \to \Cal G_\Cal U$ such that $\phi_0$ is the identity on the space of objects $T_\Cal U$ and the kernel of $\phi$ is isomorphic to  $\Cal T_\Cal U \times  \Gamma$. We say that the extension $\phi$ is topologically split if there is a continuous map $\sigma: \Cal G_\Cal U \to \tilde \Cal G$ such that $\phi(\sigma(g)) = g$ and $\sigma(1_x) = 1_x$ for every  $x \in T_\Cal U$.

 Such an extension is completely determined by a homomorphism from $\Cal G$ to the $2$-group associated to the crossed module $\mu$. The equivalence classes of such extensions is in bijection with the homotopy classes of maps from the classifying space  $B\Cal G$ of $\Cal G$  to the classifying space $B(\mu)$ of $\mu$, at least if $\Cal G$ satisfies the condition 1) in IV.2. See the Erratum below.

\subheading{3. The case where $\Gamma$ is a dense subgroup of a Lie group $G$} Let $G$ be a simply connected Lie group $G$ and  let $\Gamma$ be a dense subgoup of $G$ endoved with the discrete topology. Let $\Cal G$ be the action groupoid $\Gamma \ltimes G$, where $\Gamma$ acts on $G$ by left translations. Let $Aut(G,\Gamma)$ be the group of automorphism of $G$ preserving the subgroup $\Gamma$. The group of pointed selfequivalences of $\Cal G$ is isomorphic to the Lie group $G \rtimes Aut (G,\Gamma)$, where $Aut(G,\Gamma)$ is considered as a discrete group acting in the obvious way on the Lie group $G$. To $(g,\psi) \in G \rtimes Aut(G,\Gamma)$ corresponds the pointed self-equivalence of $\Gamma \ltimes G$ given by the $\psi$-equivariant homeomorphism $G \to G$ sending $h$ to $\psi(h)g^{-1}$.

The group $Aut(G,\Gamma)$ depends of the arithmetic properties of $\Gamma$ as a subgroup of $G$. For instance, if $G = \Bbb R$ and if $\Gamma$ is a subgroup generated by two elements whose ratio is an irrational number $\alpha$, then $Aut(G,\Gamma)$ is the group acting on $\Bbb R$ by the multiplication by $\mp 1$ when $\alpha$ is transcendental, otherwise  by the units of the ring of integers of the number field $\Bbb Q(\alpha)$.

The crossed module associated to the group of selfequivalence of $\Gamma \ltimes G$ is 
$$\mu: \Gamma \to G \rtimes Aut(G,\Gamma), $$
defined by $\mu(\gamma) =(\gamma^{-1}, Ad(\gamma))$. 

The corresponding extension problem has been studied in the thesis of Ana Maria da Silva (\cite{11}). Let $W$ be a paracompact differentiable manifold considered as a trivial \'etale groupoid. An extension of $W$ by $\Gamma \ltimes G$ is given by an open cover $\Cal U$ of $W$ and a homomorphism $\phi$ from    an etale groupoid $\tilde \Cal G$ to $W_\Cal U$. The  space of units of $\tilde \Cal G$ is assumed to be  isomorphic to $W_\Cal U \times G$  and the  kernel of $\phi$ is the action groupoid given by $\Gamma$ acting on $T_\Cal U \times G$ by left translations on the factor $G$ and trivially on the first factor.

This problem is motivated by Molino's structure theorem of Riemannian foliations on a complete Riemannian manifold $M$ which are transversally complete (see \cite 9). In that case Molino showed that the closure of the leaves are fibers of a fibration of $M$ with base space a smooth manifold $W$. The transverse holonomy groupoid of the foliation restricted to a fiber is equivalent to an action groupoid $\Gamma \ltimes G$. The transverse holonomy groupoid of the foliation is equivalent to an extension of $W$ by $\Gamma \ltimes G$ in the above sense.

Ana Maria da Silva showed in her thesis that the set of  equivalence classes of extensions of $W$ by $\Gamma \ltimes G$ are in bijection with the set of homotopy classes of maps from $W$ to a topological space $B_{(G,\Gamma)}$. This space  should be the classifying space of the crossed module $\mu$ (or equivalently the geometric realization of the nerve of the $2$-group associated to the topological crossed module $\mu$).

\subheading{4. Erratum to \cite 6} The problem of the classification of extensions of an \'etale groupoid by a discrete group was briefly discussed at the end of my expository paper \cite 6. 
We take this opportunity to 
list a few  corrections of the last pages of this paper.

\medskip

\noindent p.97  line -3   read  $\sigma:\Cal G_\Cal U \to \tilde \Cal G_\Cal  U$ such that $\phi \circ \sigma$ is the identity of 
$\Cal G_\Cal U$ and $\sigma(1_x) = 1_x$.

\noindent p.98  lines 
 3-8 should be replaced by

 The kernel $T_\Cal U \times N$ of $\phi$ is a $\tilde \Cal G_\Cal U$-sheaf of groups. For $\tilde g \in \tilde \Cal G_\Cal U$, with source $x$ and target $y$, its action is given by the relation $$\tilde g(x,n) = (y,\tilde \psi(\tilde g)(n))\tilde g,$$ where $\tilde \psi$ is a continuous homomorphism from $\tilde \Cal G_\Cal U$ to $Aut(N)$.
   Passing to the quotient we get a homomorphism $\psi : \Cal G_\Cal U \to Out(N)$ which determines a $Out(N)$ -principal bundle over $T$ with a compatible $\Cal G$-action, in other words a (generalized) morphism $\Psi:\Cal G \to Out(N)$. Let  $\Cal C$ be the associated locally constant  $\Cal G$-sheaf over $T$ with stalk $C$ (the center of $N$)  associated to the action of $Out(N)$ on $C$.

   \medskip
   Theorem 8.2 should be stated as follows:
   Let $\Psi: \Cal G \to Out(N)$ be a morphism and let $ \Cal C$ be the associated $\Cal G$- sheaf of abelian groups with stalk the center $C$ of $N$. There is a topologically split  extension of $\Cal G$ by $N$ with associated morphism $\Psi$ iff an obstruction in $\check H^3(\Cal G, \Cal C)$ vamishes. If this is the case, the set of equivalence classes of locally topologically split extensions of $\Cal G$ by $N$ with associated morphism $\Psi$ is in bijection with the set $\check H^2(\Cal G,\Cal C)$. 
 
 \medskip
\noindent  p.99  last line before the references 
 
" if and only if this homomorphism has a discrete image in $G$."

\bye

\ref \no
\by
\paper
\jour
\vol
\yr
\pages 
\endref

\ref \no
\by
\paper
\jour
\vol
\yr
\pages 
\endref

\ref \no
\by
\paper
\jour
\vol
\yr
\pages 
\endref

\ref \no
\by
\paper
\jour
\vol
\yr
\pages 
\endref

\bye